\documentclass[submission, Phys]{SciPost}

\usepackage{amssymb,amsthm,journals,microtype,mymath}
\usepackage[varg]{txfonts}

\newtheorem{theorem}{Theorem}

\usepackage{mymath}

\begin{document}

\begin{center}
  {\Large\textbf{Asymptotic Expansion of Laplace-Fourier-Type Integrals}}
\end{center}

\begin{center}
  Sara Konrad$^*$, Matthias Bartelmann
\end{center}

\begin{center}
  Institute for Theoretical Physics, Heidelberg University, Germany \\
  * sara.konrad@stud.uni-heidelberg.de
\end{center}

\begin{center}
  \today
\end{center}


\section*{Abstract}
{\bf We study the asymptotic behaviour of integrals of the Laplace-Fourier type
\begin{equation*}
  P(k) = \int_\Omega\E^{-|k|^sf(x)}\E^{\I kx}\D x\;,
\end{equation*} 
with $k\in\mathbb{R}^d$ in $d\ge1$ dimensions, with $\Omega\subset\mathbb{R}^d$ and sufficiently well-behaved functions $f:\Omega\to\mathbb{R}$. Our main result is
\begin{equation*}
  P(k) \sim \frac{\E^{-|k|^sf(0)}}{|k|^{sd/2}}\sqrt{\frac{(2\pi)^d}{\det A}}
  \exp\left(-\frac{k^\top A^{-1}k}{2|k|^s}\right)
\end{equation*}
for $|k|\to\infty$, where $A$ is the Hessian matrix of the function $f$ at its critical point, assumed to be at $x_0 = 0$. In one dimension, the Hessian is replaced by the second derivative, $A = f''(0)$. We also show that the integration domain $\Omega$ can be extended to $\mathbb{R}^d$ without changing the asymptotic behaviour.}

\vspace{10pt}
\noindent\rule{\textwidth}{1pt}
\tableofcontents\thispagestyle{fancy}
\noindent\rule{\textwidth}{1pt}
\vspace{10pt}

\newcommand{\MB}[1]{{\color{blue}\emph{#1}}}
\newcommand{\SK}[1]{{\color{orange}\emph{#1}}}

\section{Introduction}
\label{sec:1}

Integrals over rapidly oscillating integrands occur frequently in physics and are often difficult to handle numerically. Their asymptotic behaviour for large parameter values (e.g.\ large wave numbers) is interesting in view of estimates, studies of limiting behaviour, and testing numerical solutions. In this paper, we study the asymptotic behaviour of integrals of the form
\begin{equation}
  P(k) = \int_\Omega\E^{-k^sf(x)}\E^{\I kx}\D x\;,\quad k > 0\;,
\label{eq:1}
\end{equation}
under fairly general assumptions for $k\to\infty$, first in one, then in $d$ dimensions. Specifically, we shall assume that the integral in (\ref{eq:1}) converges absolutely and that $s\ge\alpha$, with $\alpha$ to be defined below. Our main results, derived in Sect.~\ref{sec:2} for $d = 1$ and in Sect.~\ref{sec:3} for $d>1$, will be:
\begin{enumerate}
  \item In one dimension,
  \begin{equation}
    P(k) \sim \frac{\E^{-k^sf(0)}}{k^{s/2}}
    \sqrt{\frac{2\pi}{f''(0)}}
    \exp\left(-\frac{k^{2-s}}{2f''(0)}\right)
    \quad\mbox{for}\quad k\to\infty\;,
  \label{eq:2}
  \end{equation}
  where $x_0 = 0$ is a critical point of $f$.
  \item In $d>1$ dimensions,
  \begin{equation}
    P(k) \sim \frac{\E^{-|k|^sf(0)}}{|k|^{sd/2}}\sqrt{\frac{(2\pi)^d}{\det A}}
    \exp\left(-\frac{k^\top A^{-1}k}{2|k|^s}\right)
    \quad\mbox{for}\quad k\to\infty\;,
  \label{eq:3}
  \end{equation}
  where $A$ is the Hessian matrix of $f$ at its critical point, assumed to be in $x_0 = 0$.
\end{enumerate}

With this result, we extend the existing literature in two ways: we derive the asymptotic behaviour of Laplace-Fourier type integrals in one dimension and extend them to an arbitrary number $d$ of dimensions. In one dimension, detailed studies exist with a different focus \cite{berg1962, pederson1965, schell1978}, and to our knowledge, the generalization to $d$ dimensions is new.

\section{One-dimensional case}
\label{sec:2}

\subsection{Erdélyi's theorem for Laplace integrals}

In one dimension, we begin with Erdélyi's theorem \cite{doetsch1955handbuch, UBHD-65249153, 1961ArRMA...7....1E} for Laplace integrals, which states:

\begin{theorem}[Erdélyi]
\label{the:1}
  Let $I(\lambda)$ be an integral of the form
  \begin{equation}
    I(\lambda) = \int_a^b\E^{-\lambda f(x)}g(x)\,\D x\;,
  \label{eq:4}
  \end{equation}
  where $f(x)$ is a real function of the real variable $x$, while $g(x)$ may be real or complex. Then, if
  \begin{enumerate}
    \item $f(x)>f(a)$ for $x\in(a,b)$ and
    \begin{equation}
      \inf_{[a+\delta,b)}f(x)-f(a)>0
    \label{eq:5}
    \end{equation}
    for $\delta>0$;
    \item $f'(x)$ and $g(x)$ are continuous in a neighbourhood of $a$, except possibly at $a$;
    \item $f$ and $g$ admit asymptotic expansions
    \begin{align}
      f(x) &\sim f(a)+\sum_{k=0}^\infty a_k(x-a)^{k+\alpha}\;,\nonumber\\
      g(x) &\sim \sum_{k=0}^\infty b_k(x-a)^{k+\beta-1}\;,
    \label{eq:6}
    \end{align}
    $f$ can be term-wise differentiated,
    \begin{equation}
      f'(x) \sim \sum_{k=0}^\infty a_k(k+\alpha)(x-a)^{k+\alpha-1}
    \label{eq:7}
    \end{equation}
    for $x\to a^+$, where $\alpha>0$ and $\mathrm{Re}\,\beta>0$; and
    \item $I(\lambda)$ converges absolutely for sufficiently large $\lambda$; then
  \end{enumerate}
  the integral $I(\lambda)$ has the asymptotic expansion
  \begin{equation}
    I(\lambda) \sim \E^{-\lambda f(a)}\sum_{n=0}^\infty
    \frac{\Gamma(\nu)\,c_n}{\lambda^\nu}
  \label{eq:8}
  \end{equation}
  for $\lambda\to\infty$, where $\nu := (n+\beta)/\alpha$. The coefficients $c_n$ can be expressed by $a_n$ and $b_n$ as
  \begin{equation}
    c_n = \frac{1}{\alpha a_0^\nu}\sum_{m=0}^n
    \frac{b_{n-m}}{m!}\,d_{m,n}\quad\mbox{with}\quad
    d_{m,n} = \lim_{x\to0}\frac{\D^m}{\D x^m}\left(
      1+\sum_{j=1}^\infty\frac{a_j}{a_0}\,x^j
    \right)^{-\nu}\;.
  \label{eq:9}
  \end{equation}
\end{theorem}

Beginning with this result, we insert the derivatives $d_{m,n}$ into the coefficients $c_n$ from (\ref{eq:9}) and write the asymptotic expansion $I(\lambda)$ from (\ref{eq:8}) for some positive integer $M$ that is to be constrained later, as
\begin{align}
  I(\lambda) &\sim \frac{\E^{-\lambda f(a)}}{\alpha}
  \sum_{n=0}^{\infty} \sum_{m=0}^n
  \frac{\Gamma(\nu)}{(a_0\lambda)^\nu}\,\frac{bac_{n-m}}{m!}\,d_{m,n}
  \nonumber\\ &=
  \frac{\E^{-\lambda f(a)}}{\alpha} \left[
  \sum_{m=0}^{M-1}\sum_{n=m}^{\infty}
  \frac{\Gamma(\nu)}{(a_0\lambda)^\nu}\,\frac{b_{n-m}}{m!}\,d_{m,n} + \phi_M(\lambda) \right]\;,
\label{eq:10}
\end{align}
where we have assumed that the inner sum over $n$ in the second line exists for all $0\le m\le M$. To abbreviate notation, we have defined
\begin{equation}
    \phi_M(\lambda) := \sum_{n=M}^{\infty}\sum_{m=M}^n \frac{\Gamma(\nu)}{(a_0\lambda)^\nu}\,\frac{b_{n-m}}{m!}\,d_{m,n}\;.
\label{eq:11}
\end{equation}
and introduce the symbol
\begin{equation}
  I_m(\lambda) := \sum_{n=m}^\infty
  \frac{\Gamma(\nu)}{(a_0\lambda)^\nu}\,\frac{b_{n-m}}{m!}\,d_{m,n}
\label{eq:12}
\end{equation} 
for the sums over $n$ remaining in (\ref{eq:10}). If the outer sum over $m$ in (\ref{eq:10}) exists for $M\to\infty$, can write the asymptotic expansion $I(\lambda)$ as
\begin{equation}
  I(\lambda) \sim \frac{\E^{-\lambda f(a)}}{\alpha}
    \sum_{m=0}^\infty I_m(\lambda)\;.
\label{eq:13}
\end{equation}
Note that this is indeed an asymptotic expansion, since the symbols $I_m(\lambda)$ are an asymptotic sequence as $\lambda  \rightarrow \infty$ in $\mathbb{R}^+$, denoted by $\{I_m\}$, because for all $m \ge 0$,
\begin{equation}
  I_{m+1}(\lambda) = o\left(I_m(\lambda)\right)
  \quad\text{as}\quad \lambda \rightarrow \infty\;;
\label{eq:14}
\end{equation}
see \cite[p.~10, definition 2]{UBHD-65467105}. The sum over $n$ in the functions $I_m(\lambda)$ defined in (\ref{eq:12}) can be carried out for some choices of $\alpha$ and $\beta$, for which we shall give examples below.

\subsection{Laplace-Fourier integrals in one dimension}

We now wish to study the asymptotic behaviour of the Fourier integral
\begin{equation}
  I(\lambda,k) = \int_\Omega \E^{-\lambda f(x)}x^{\,\beta-1}\E^{\I k x}\,\D x
\label{eq:15}
\end{equation}
for $\lambda\to\infty$, with $\mathrm{Re}\,\beta>0$ and the fixed parameter $k > 0$. For doing so, we apply Erdélyi's theorem, setting
\begin{equation}
  g(x) = x^{\beta-1}\E^{\I k x}\;.
\label{eq:16}
\end{equation} 
Moreover, without loss of generality, we shift the origin such that $a = 0$ and take the integral $I(\lambda,k)$ over the open domains $\Omega=(0,b)$ and assume that $f(x)$ satisfies the conditions of Erdélyi's theorem with $a = 0$. The coefficients $b_n$ defined in (\ref{eq:6}) are given by
\begin{equation}
  b_n(k) = \frac{(\I k)^n}{n!}\;,
\label{eq:17}
\end{equation}
and the symbols $I_m$ defined in (\ref{eq:12}) are given by
\begin{equation}
  I_m(\lambda,k) = \frac{1}{m! (a_0 \lambda)^{(m+\beta)/\alpha}} \sum_{n=0}^{\infty} \frac{\Gamma\left( \nu + \frac{m}{\alpha}\right) }{n!}\left( \frac{\I k}{(a_0 \lambda)^{1/\alpha}} \right)^n d_{m,n+m}\;.
\label{eq:18}
\end{equation}

We now show that, for fixed $\lambda, k$, the symbols $I_m$ as in (\ref{eq:18}) converge for each finite $m$ when $\alpha > 1$. For $\alpha = 1$, the symbols converge if $\left| \frac{k}{a_0 \lambda}\right| < 1$.

\begin{proof}
For $m > 0$, the derivatives $d_{m,n}$ as defined in (\ref{eq:9}) can be written as
\begin{equation}
  d_{m,n} = \sum_{p=1}^m \frac{\Gamma\left(\nu+p\right)}{\Gamma\left(\nu\right)}\, D_{m,p}\left(a_0,\dotsm,a_{m-p+1}\right)\;,
\label{eq:19}
\end{equation}
where $D_{m,p}$ are finite coefficients that are independent of $n$. To see this, note that any of the $m$ derivatives applied in (\ref{eq:9}) acts on the parenthesis raised to the power of $-\nu$, or on the sum in the parenthesis. The former derivatives create the Gamma functions in front, the latter the coefficients $D_{m,p}$. As the limit $x\to0$ is taken, the expression in the parenthesis becomes unity, thus raising it to any power becomes independent of the exponent.\footnote{A systematic way to compute these coefficients in terms of partial ordinary Bell polynomials is presented in \cite{2013ConAp..38..471N}.} Thus, for each $m$ there exists a positive constant $K_m > 0$ such that the following estimate holds
\begin{equation}
  \left| d_{m,n}\right| \le K_m \sum_{p=1}^m \frac{\Gamma\left(\nu+p\right)}{\Gamma\left(\nu\right)} \le m \cdot K_m  \frac{\Gamma\left(\nu+m\right)}{\Gamma\left(\nu\right)}\;.
\label{eq:20}
\end{equation}
Using (\ref{eq:17}), (\ref{eq:18}) and (\ref{eq:20}), the symbols $I_m$ are then bound by
\begin{equation}
  \left|I_m(\lambda,k)\right| \le \frac{K_m}{(m-1)!(a_0 \lambda)^{\frac{m+\beta}{\alpha}}} \sum_{n=0}^\infty \frac{\Gamma \left( \nu + m + \frac{m}{\alpha} \right)}{\Gamma\left(n+1\right)}\left|\frac{k}{(a_0 \lambda)^{1/\alpha}}\right|^n\;.
\label{eq:21}
\end{equation}
For $\alpha > 1$, the sum in eq. (\ref{eq:21}) is absolutely convergent, which can be seen by d'Alembert's ratio test. We consider the limit of the ratio of two consecutive terms
\begin{align*}
    & \lim_{n\rightarrow \infty} \left| \frac{\Gamma \left( \frac{n+1+\beta+(1+\alpha)m}{\alpha} \right)}{\Gamma(n+2)} \left| \frac{k}{(a_0 \lambda)^{\frac{1}{\alpha}}} \right|^{n+1} \cdot \frac{\Gamma \left(n+1\right)}{\Gamma \left( \frac{n+\beta+(1+\alpha)m}{r+1} \right)} \left| \frac{(a_0 \lambda)^{\frac{1}{\alpha}}}{k} \right|^{n} \right|
    \\
    &= \lim_{n \rightarrow \infty} \left| \frac{k}{(a_0 \lambda)^{\frac{1}{\alpha}}} \right| \frac{1}{n+1} \frac{\Gamma \left( \frac{n+1+(1+\alpha)m+\beta-\alpha}{\alpha}+1 \right)}{\Gamma \left( \frac{n+1+(1+\alpha)m+\beta-\alpha}{\alpha}+1- \frac{1}{\alpha} \right)}\\
    & \le \frac{k}{(a_0 \lambda)^{\frac{1}{\alpha}}} \lim_{n \rightarrow \infty} 
    \frac{1}{n+1} \left( \frac{n+1+(1+\alpha)m+\beta}{\alpha} \right)^{\frac{1}{\alpha}} \\
    &=  \frac{k}{(a_0 \lambda)^{\frac{1}{\alpha}}} \left(\frac{1}{\alpha}\right)^{\frac{1}{\alpha}} 
    \lim_{n \to 0}n^{-1+\frac{1}{\alpha}} 
    = 0\;,
\end{align*}
where we made use of Gautschi's inequality $x^{1-s} < \frac{\Gamma (x+1)}{\Gamma(x+s)} < (x+1)^{1-s}$, for $0 < s < 1$ and $x > 0$ \cite{10.5555/1830479}. Thus, the sums (\ref{eq:21}) are finite for $m > 0$ and $\alpha > 1$. Since $d_{0,n}=1$ for all $n$, an analogous proof shows $|I_0(\lambda,k)| < \infty$.
For $\alpha = 1$, the inequality (\ref{eq:21}) becomes
\begin{align}
     \left|I_m(\lambda,k)\right| 
     &\le \frac{K_m}{(m-1)!(a_0 \lambda)^{m+\beta}} \sum_{n=0}^\infty \frac{\Gamma \left( \nu + 2m \right)}{n!}\left|\frac{k}{a_0 \lambda}\right|^n \nonumber\\
     &= \frac{K_m}{(m-1)!(a_0 \lambda)^{m+\beta}} \Gamma(2m+ \beta) {}_1F_0\left(2m+\beta;; \left| \frac{k}{a_0 \lambda} \right|\right)\;,
\label{eq:22}
\end{align}
where the last step is valid for $\left| \frac{k}{a_0 \lambda} \right| < 1$ and ${}_pF_q$ denotes the generalized hypergeometric function. The case $m = 0$ is easier because $d_{0,n} = 1$ for all $n$, but otherwise proceeds analogously.
\end{proof}

Now that we have shown, that (\ref{eq:13}) together with the symbols $I_m(\lambda, k)$ is a valid asymptotic expansion for Laplace-Fourier integrals $I(\lambda, k)$ as $\lambda \to \infty$, we compare our expansion (\ref{eq:13}) to the expansion (\ref{eq:9}) from Erdélyi's theorem.
Since the symbols $I_m$ originate from a resummation of Erdélyi's series, they contain the coefficients $b_n$ to infinite order, thus the total information about the phase function $\E^{\I kx}$.
Considering now Laplace-Fourier integrals with related exponential factors  $\lambda = k^s$ implies $P(k) = I(k^s,k)$. The symbols $I_m$ are then specialized to
\begin{equation}
  I_m(k^s,k) = \frac{1}{m! (a_0 k^s)^{(m+\beta)/\alpha}} \sum_{n=0}^{\infty} \Gamma\left( \nu + \frac{m}{\alpha}\right) \left( \frac{\I k^{1-s/\alpha}}{a_0^{1/\alpha}} \right)^n d_{m,n+m}\;.
\label{eq:23}
\end{equation}

We find for $s \ge \alpha$
\begin{equation}
  I_m(k^s, k) = \mathcal{O}\left(k^{-\frac{sm+s\beta}{\alpha}}\right)\;, \quad \text{as } k \to \infty\;,
\label{eq:24}
\end{equation}
showing that the symbols $I_m(k^s,k)$ denote an asymptotic sequence as $k \rightarrow \infty$ in $\mathbb{R}^+$. Comparing this sequence to the terms in the expansion (\ref{eq:9}) from Erdélyi's theorem,
\begin{equation}
  \frac{\Gamma(\nu)c_n}{k^{s\nu}} =
  \mathcal{O}\left(k^{-\frac{sn+s\beta}{\alpha}+n}\right)
  \quad\text{as}\quad k \to \infty\;,
\label{eq:25}
\end{equation}
shows that the terms in our asymptotic expansion (\ref{eq:13}) fall of much faster when $s > \alpha$. For $s = \alpha$, the sequence from Erdélyi's theorem even fails to be asymptotic.

\subsection{Examples}

For convenience, we introduce
\begin{equation}
  x_\alpha := \frac{4}{a_0\lambda}\left(\frac{k}{4}\right)^\alpha\;.
\label{eq:26}  
\end{equation}
Note that for $\lambda = k^{s}$ and $s=\alpha$, $x_\alpha$ is a positive constant, while for $s > \alpha$, $x_\alpha > 0$ approaches $0$ in the limit $k \to \infty$.

\subsubsection{One-sided domain $\Omega = (0,b)$}

On the one-sided domain $\Omega$, we can then directly apply Erdélyi's theorem by inserting the $b_n$ from (\ref{eq:17}) into (\ref{eq:12}). Since $d_{0,n} = 1$ for all $n$, we then have
\begin{equation}
  I_0(\lambda, k) = \frac{1}{(a_0 \lambda)^{\beta/\alpha}} 
  \sum_{n=0}^\infty
  \frac{\Gamma(\nu)}{n!}\left(\frac{\I k}{(a_0 \lambda)^{1/\alpha}}\right) 
  = 
  \frac{1}{(a_0 \lambda)^{\beta/\alpha}} \
  {}_1\Psi_0 \left[ 
  \begin{smallmatrix} (\beta/\alpha, 1/\alpha) \\ -
  \end{smallmatrix}; \I x_1 \right]\;
\label{eq:27}
\end{equation}
for $m = 0$, where ${}_p\Psi_q$ denotes the Fox-Wright function. For this sum to converge, we must require $\alpha\ge1$ and thus $\nu\le n+\beta$. For $\alpha = 1$ and arbitrary $\beta>0$, we arrive at
\begin{equation}
  I_0(\lambda, k) =
  \Gamma(\beta)\left(\frac{x_1}{k\,(1-\I x_1)}\right)^{\,\beta}\;.
\label{eq:28}
\end{equation}

\subsubsection{Two-sided domain $\Omega = (-b_1, b_2)$}

On the two-sided domain $\Omega$, we assume that Erdélyi's theorem applies to $I(\lambda, k)$ separately on the intervals $(0,b_1)$ with $x\to-x$ and $(0, b_2)$. Furthermore, we restrict $\alpha$ to being an even integer because otherwise $f(x)$ could not have a maximum in the origin. In this case, it is advantageous to note that $b_n \to (-1)^{\,\beta-1+n}b_n$ and $d_{m,n}\to(-1)^md_{m,n}$ for $x\to-x$.

Then, after adding the results for the intervals $(-b_1, 0)$ and $(0, b_2)$, only those summands in $I_m(\lambda,k)$ from (\ref{eq:12}) remain for integer $\beta$ which have even $n$ if $\beta = 2l-1$ is odd, or which have odd $n$ if $\beta = 2l$ is even, where $l\in \{1,2,\dots\}$. We use this to write
\begin{equation}
  I_m(\lambda, k) =
  2\,\begin{cases}
    \displaystyle
    \sum_{n=\lceil m/2\rceil}^\infty
      \frac{\Gamma\left(\frac{2n+\beta}{\alpha}\right)}
           {(a_0\lambda)^{(2n+\beta)/\alpha}}
      \frac{(\I k)^{2n-m}}{m!(2n-m)!}d_{m,2n} & \beta = 2l-1 \\[12pt]
    \displaystyle
    \sum_{n=\lfloor m/2\rfloor}^\infty
      \frac{\Gamma\left(\frac{2n+1+\beta}{\alpha}\right)}
          {(a_0\lambda)^{(2n+1+\beta)/\alpha}}
      \frac{(\I k)^{2n+1-m}}{m!(2n+1-m)!}d_{m,2n+1} & \beta = 2l \\
  \end{cases}\;,
\label{eq:29}
\end{equation}
with the factor of $2$ accounting for the results on both intervals.

Some special cases for the leading term $I_0(\lambda)$ may be worth noting:

\begin{enumerate}
  \item[(i)] For $\alpha = 2$ and integer $\beta$, the asymptotic expansion of $f(x)$ begins with a term quadratic in $x$, and the envelope of the Fourier factor in (\ref{eq:16}) has an integer power of $x$. In this case, we have
  \begin{equation}
    I_0(\lambda,k) =
    \begin{cases}
      \displaystyle
      2\left(\frac{4x_2}{k^2}\right)^{l+1/2}\,
        \Gamma\left(l+\frac{1}{2}\right)\,
        {}_1F_1\left(l+\frac{1}{2};\frac{1}{2};-x_2\right) & \beta = 2l+1 \\
      \displaystyle
      2\I k\left(\frac{4x_2}{k^2}\right)^{l+1/2}\,
        \Gamma\left(l+\frac{1}{2}\right)\,
        {}_1F_1\left(l+\frac{1}{2};\frac{3}{2};-x_2\right) & \beta = 2l \\
    \end{cases}\;.
  \label{eq:30}
  \end{equation} 
  
  \item[(ii)] For $\alpha = 4$ and $\beta = 1$, the asymptotic expansion of $f(x)$ begins with a term quartic in $x$, and the Fourier factor has a constant amplitude in $x$. Then, the result is
  \begin{equation}
    I_0(\lambda,k) = 2\sqrt{2}\,\frac{x_4^{1/4}}{k}
    \left[
      2\Gamma\left(\frac{1}{4}\right)
      {}_0F_2\left(;\frac{1}{2},\frac{3}{4};\frac{x_4}{4}\right)-
      8x_4^{1/2}\Gamma\left(\frac{3}{4}\right)
      {}_0F_2\left(;\frac{5}{4},\frac{3}{2};\frac{x_4}{4}\right)
    \right]\;.
  \label{eq:31}
  \end{equation}
  In the expressions (\ref{eq:30}) and (\ref{eq:31}), ${}_pF_q$ is the generalised hypergeometic function.
\end{enumerate}

In view of applications in physics, we investigate in more detail situations in which $f(x)$ has an asymptotic expansion beginning with a term quadratic in $x$, and the amplitude of the Fourier factor is constant. Then, $\alpha = 2$ and $\beta = 1$. In this case, we return directly to the definition (\ref{eq:12}) for the summands $I_m(\lambda,k)$ in the asymptotic expansion (\ref{eq:13}) of the integral $I(\lambda,k)$. Since $\beta$ is odd now, we can restricting the sum over $n$ to even $n$, and since the integration domain is two-sided, we need to double the result. Then, from (\ref{eq:12}),
\begin{equation}
  I_0(\lambda,k) = 2\sum_{n=0}^\infty
  \frac{\Gamma\left(n+\frac{1}{2}\right)}
       {(a_0\lambda)^{n+1/2}}\frac{(\I k)^{2n}}{(2n)!} =
  2\sqrt{\frac{\pi}{a_0\lambda}}
  \sum_{n=0}^\infty \frac{1}{n!} \left( - \frac{k^2}{4a_0\lambda} \right)^n =
  2\sqrt{\frac{\pi}{a_0\lambda}}\,\E^{-x_2}\;,
\label{eq:32}
\end{equation}
where we have used
\begin{equation}
  \Gamma\left(n+\frac{1}{2}\right) = \frac{(2n)!}{n!4^n}\sqrt{\pi}
\label{eq:33}
\end{equation}
in the second step. With
\begin{equation}
  d_{1,2n} = -\frac{2n+1}{2}\frac{a_1}{a_0}\quad\mbox{and}\quad
  d_{2,2n} = (2n+1)\left[
    \frac{2n+3}{4}\frac{a_1^2}{a_0^2}-\frac{a_2}{a_0}
  \right]\;,
\label{eq:34}
\end{equation}
we further find
\begin{equation}
  \frac{I_1(\lambda,k)}{I_0(\lambda,k)} = \frac{\I A}{k} \quad\mbox{and}\quad
  \frac{I_2(\lambda,k)}{I_0(\lambda,k)} = \frac{B}{k^2}
\label{eq:35}
\end{equation}
with
\begin{align}
  A &:= \frac{a_1}{a_0}\,x_2\left(2x_2-3\right)\quad\mbox{and}\nonumber\\
  B &:= -\frac{a_2}{a_0}\,x_2\left(4x_2^2-12x_2+3\right)-
  \frac{a_1^2}{4a_0^2}\,x_2\left(8x_2^3-60x_2^2+90x_2-15\right)\;.
\label{eq:36}
\end{align}

We finally specialise this result to integrals of the form
\begin{equation}
  P(k) = \int_\Omega\E^{-k^sf(x)}\E^{\I k x}\,\D x 
\label{eq:37}
\end{equation}
with $s \ge 2$ and $f$ satisfying the previously stated conditions. Substituting $\lambda\to k^s$ in (\ref{eq:32}) then immediately leads to the asymptotic series
\begin{equation}
  P(k) \sim \frac{\E^{-k^sf(0)}}{k^{s/2}}\sqrt{\frac{\pi}{a_0}}\,
  \E^{-x_2}\left[
    1+\I\,\frac{A}{k}+\frac{B}{k^2}+\mathcal{O}\left(k^{-3}\right)
  \right]
\label{eq:38}
\end{equation}
for $k\to\infty$, with $x_2$ specialized to $x_2 = k^{2-s}/(4a_0)$. Note that the terms $A$ and $B$ in (\ref{eq:34}) are both of order $\mathcal{O}(k^{s-2})$ here. The result (\ref{eq:2}) follows from (\ref{eq:34}) by noting that $a_0 = f''(0)/2$.

\section{$d$-dimensional case}
\label{sec:3}

Before we can now turn to the case of an arbitrary number $d$ of dimensions, we need to prepare some concepts and notation. Here and below, $D$ denotes the derivative of a function with respect to its arguments.

\subsection{Preliminary remarks}

Let $\Omega\subset\mathbb{R}^d$ be an open subset of $\mathbb{R}^d$ and $x_0\in\Omega$ a non-degenerate critical point of a smooth function $f:\Omega\to\mathbb{R}$ on $\Omega$ with $f\in C^\infty$. Then, by Morse's lemma, neighbourhoods $U$, $V$ of the points $y = 0$ and $x = x_0$ and a diffeomorphism $h: U\to V$ exist such that
\begin{equation}
  \left(f\circ h\right)(y) = f(x_0)+\frac{1}{2}y^\top Qy\;,\quad
  Q = \mathrm{diag}\left(\mu_1,\ldots,\mu_d\right)\;,
\label{eq:39}
\end{equation}
where the Jacobian $H := Dh$ of $h$ has unit determinant at $y = 0$, $\det H_0= 1$ with $H_0 := H(0)$. Then, by the chain rule and $(Df)(x_0) = 0$, $H_0$ diagonalises the Hessian matrix
\begin{equation}
  A := \left(D^2f\right)\Bigr\vert_{x=x_0}
\label{eq:40}
\end{equation}
of $f$ in $x_0$ such that
\begin{equation}
  Q = H^\top_0AH_0\;.
\label{eq:41}
\end{equation}

With the same chart $h:U\to V$, the Jacobian $H=Dh$ and a smooth function $g:\Omega\to\mathbb{R}$ with $g\in C^\infty$, we define the function $G:U\to\mathbb{R}$ by $G(y) = (g\circ h)(y)\det H(y)$. We further introduce multi-indices $\alpha=(\alpha_1, \alpha_2, \ldots, \alpha_d)$ and agree on the notation
\begin{align}
  |\alpha| &:= \alpha_1+\alpha_2+\ldots+\alpha_d\;,\nonumber\\
  \alpha!  &:= \alpha_1!\alpha_2!\cdots\alpha_d!\;,\nonumber\\
  \Gamma(\alpha) &:=
  \Gamma(\alpha_1)\Gamma(\alpha_2)\cdots\Gamma(\alpha_d)\;,\nonumber\\
  \mu^\alpha &:= \mu_1^{\alpha_1}\mu_2^{\alpha_2}\cdots\mu_d^{\alpha_d}\;,
\label{eq:42}
\end{align}
where $\mu = (\mu_1, \mu_2, \ldots, \mu_d)$ is a $d$ dimensional real valued vector.
We finally define the symbol
\begin{equation}
  \delta(\alpha) = \begin{cases}
                     1 & \mbox{all $\alpha_j$ even} \\
                     0 & \mbox{else}
                   \end{cases}
\label{eq:43}
\end{equation}
and the derivative operator
\begin{equation}
  D^\alpha G(0) = \frac
  {\partial^{|\alpha|}}{\partial^{\alpha_1}y_1\ldots\partial^{\alpha_d}y_d}
  G(y)\Biggr\vert_{y = 0}\;.
\label{eq:44}
\end{equation}
With these preparations, we can now continue with the following theorem for the asymptotic behaviour of a class of integrals on $\Omega$ \cite{UBHD-65467105}.

\begin{theorem}
\label{the:2}
  Let $J(\lambda)$ be the integral
  \begin{equation}
    J(\lambda) = \int_\Omega\E^{-\lambda f(x)}g(x)\,\D x
  \label{eq:45}
  \end{equation}
  with $x\in\mathbb{R}^d$ and $\Omega\subset\mathbb{R}^d$, whose asymptotic behaviour we wish to determine for $\lambda\to\infty$. We assume that
  \begin{enumerate}
    \item  $f, g\in C^\infty$ on $\Omega$;
    \item $J(\lambda)$ converges absolutely for sufficiently large $\lambda$;
    \item $f$ has a minimum at, and only at, $x_0\in\Omega$ such that
    \begin{equation}
      \rho(\varepsilon) := \inf_{\Omega\setminus B_\varepsilon(x_0)}f(x)-f(x_0) > 0
    \label{eq:46}
    \end{equation}
    for all $\varepsilon>0$, where $B_\varepsilon(x_0)$ is the open ball with radius $\varepsilon$ around $x_0$; and that
    \item the Hessian matrix $A$ of $f$ in $x_0$ is positive definite.
  \end{enumerate}
  Then, the integral $J(\lambda)$ has the asymptotic expansion
  \begin{equation}
    J(\lambda) \sim \E^{-\lambda f(x_0)}\sum_{n=0}^\infty
    \frac{c_n}{\lambda^{d/2+n}}
  \label{eq:47}
  \end{equation}
  for $\lambda\to\infty$, with
  \begin{equation}
    c_n = \sum_{|\alpha|=2n}\delta(\alpha)
    \left(\frac{2}{\mu}\right)^{(\alpha+1)/2}
    \Gamma\left(\frac{\alpha+1}{2}\right)\frac{D^\alpha G(0)}{\alpha!}\;.
  \label{eq:48}
  \end{equation}
  The $\mu_i$, $1\le i\le d$, are the (real and positive) eigenvalues of $A$.
\end{theorem}

\subsection{Laplace-Fourier integrals in $d$ dimensions}

We apply the result (\ref{eq:47}) with the coefficients (\ref{eq:48}) to Laplace-Fourier-type integrals,
\begin{equation}
  J(\lambda, k) = \int_\Omega\E^{-\lambda f(x)}\E^{\I kx}\D x\;,
\label{eq:49}
\end{equation}
whose asymptotic behaviour we wish to determine for $\lambda\to\infty$, where $|k|$ is assumed to be a large parameter at the same time. Without loss of generality, we set $x_0 = 0$ here. For integrals of this type, the function $G$ defined above is
\begin{equation}
  G(y) = \E^{\I kh(y)}\det H(y)\;.
\label{eq:50}
\end{equation}
Since only the derivatives of $G$ acting on the exponential factor increase the powers of $k$, and since additionally $G(0) = 1$, we can state
\begin{equation}
  D^\alpha G(0) \sim
  \prod_{j=1}^d\left[\I\left(H^\top_0k\right)_j\right]^{\alpha_j} =
  \I^{|\alpha|}
  \prod_{j=1}^d\left(H^\top_0k\right)_j^{\alpha_j}
  \quad\mbox{for}\quad |k|\to\infty
\label{eq:51}
\end{equation}
and move on to write
\begin{equation}
  c_n(k) \sim \sum_{|\alpha|=2n}\delta\left(\alpha\right)\left(-1\right)^n
  \prod_{j=1}^d\left(\frac{2}{\mu_j}\right)^{(\alpha_j+1)/2}
  \Gamma\left(\frac{\alpha_j+1}{2}\right)
  \frac{\left(H^\top_0k\right)_j^{\alpha_j}}{\alpha_j!}
  \quad\mbox{for}\quad |k|\to\infty\;.
\label{eq:52}
\end{equation}
Since the factor $\delta(\alpha)$ ensures that only even $\alpha_j$ can occur in the sum, it is convenient to define a multi-index $n$ with elements $n_j=\alpha_j/2$. Inserting (\ref{eq:33}) into (\ref{eq:52}) and noting that $\prod\mu_j = \det A$, we arrive at
\begin{equation}
  c_n(k) \sim \sqrt{\frac{(2\pi)^d}{\det A}}\left(-\frac{1}{2}\right)^n
  \sum_{|n|=n}\prod_{j=1}^d\frac{1}{n_j!}
  \frac{\left(H^\top_0k\right)_j^{2j}}{\mu_j^{n_j}}
  \quad\mbox{for}\quad |k|\to\infty\;.
\label{eq:53}
\end{equation}

Since (\ref{eq:41}) implies that $H_0Q^{-1}H_0^\top = A^{-1}$, the remaining sum equals
\begin{equation}
  \sum_{|n|=n}\prod_{j=1}^d\frac{1}{n_j!}
    \frac{\left(H^\top_0k\right)_j^{2j}}{\mu_j^{n_j}} =
  \frac{1}{n!}
  \left[
    \left(H^\top_0k\right)^\top Q^{-1}\left(H^\top_0k\right)
  \right]^n = \frac{\left(k^\top A^{-1}k\right)^n}{n!}
\label{eq:54}
\end{equation}
allowing us to write (\ref{eq:47}) in the form
\begin{equation}
  J(\lambda,k) \sim \E^{-\lambda f(0)}\sqrt{\frac{(2\pi)^d}{\lambda^d\det A}}
  \exp\left(-\frac{k^\top A^{-1}k}{2\lambda}\right)
  \quad\mbox{for}\quad \lambda\to\infty\;.
\label{eq:55}
\end{equation}

For integrals of the type
\begin{equation}
  P(k) = \int_\Omega\E^{-|k|^sf(x)}\E^{\I kx}\D x
\label{eq:56}
\end{equation}
with $s\ge2$, we substitute again $\lambda\to|k|^s$ in (\ref{eq:50}). We then obtain immediately from (\ref{eq:55})
\begin{equation}
  P(k) \sim \E^{-|k|^sf(0)}
  \sqrt{\frac{(2\pi)^d}{|k|^{sd}\det A}}
  \exp\left(-\frac{k^\top A^{-1}k}{2|k|^s}\right)
  \quad\mbox{for}\quad |k|\to\infty\;,
\label{eq:57}
\end{equation}
which proves the result (\ref{eq:3}).

\subsection{Extension to $\mathbb{R}^d$}

We finally wish to extend the integration domain $\Omega$ to $\mathbb{R}^d$. We begin with the conditions imposed in Theorem \ref{the:2}, but relax conditions (1) and (2) to
\begin{enumerate}
  \item[$1'$.] the Hessian $A=D^2f$ exists in $x_0 = 0$ and
  \item[$2'$.] $f\in L_2(\mathbb{R}^d)$
\end{enumerate}
and add as a further condition:
\begin{enumerate}
  \item[3b.] The infimum
  \begin{equation}
    \sigma(\varepsilon) := \inf_{\mathbb{R}^d\setminus B_\varepsilon(0)}|f(0)|-|f(x)| > 0
  \label{eq:58}
  \end{equation}
  for all $\varepsilon>0$.
\end{enumerate}
Note that conditions (3) and (3b) together imply that $f(0)<0$. Then, the integral
\begin{equation}
  P(k) = \int_{\mathbb{R}^d}\E^{-|k|^sf(x)}\E^{\I kx}\D x
\label{eq:59}
\end{equation}
has the asymptotic expansion (\ref{eq:57}).

To prove this statement, we consider the integral
\begin{align}
  \bar P(k) :&= \int_{\mathbb{R}^d}\left[\E^{-|k|^sf(x)}-1+|k|^sf(x)\right]
  \E^{\I kx}\D x \nonumber\\ &=
  P(k)-(2\pi)^N\delta_\mathrm{D}(k)+|k|^s\tilde f(k) \nonumber\\ &=
  P(k)+|k|^s\tilde f(k)
\label{eq:60}
\end{align}
for $k>0$, where $\tilde f$ denotes the Fourier transform of $f$. We further split the integration domain as
\begin{equation}
  \bar P(k) = \left(\int_\Omega+\int_{\mathbb{R}^d\setminus\Omega}\right)
  \left[
    \E^{-|k|^sf(x)}-1+|k|^sf(x)
  \right]\E^{\I kx}\D x =: \bar P_\Omega(k)+\bar P_{\bar\Omega}(k)\;.
\label{eq:61}
\end{equation}
For the asymptotic expansion of $\bar P_\Omega(k)$, we can estimate
\begin{equation}
  \bar P_\Omega(k) \le \int_\Omega\E^{-|k|^sf(x)}\E^{\I kx}\D x+\left[
    1+|k|^s|f(0)|
  \right]V_\Omega\;,
\label{eq:62}
\end{equation} 
where $V_\Omega$ is the volume of the domain $\Omega$. Since the second term on the right-hand side of (\ref{eq:62}) is exponentially suppressed compared to the asymptotic expansion (\ref{eq:57}) of the first term (recall that $f(0)<0$), we can conclude that
\begin{equation}
  \bar P_\Omega(k) \sim \E^{-|k|^sf(0)}\sqrt{\frac{(2\pi)^d}{|k|^{sd}\det A}}
  \exp\left(-\frac{k^\top A^{-1}k}{2|k|^s}\right)
  \quad\mbox{for}\quad |k|\to\infty\;.
\label{eq:63}
\end{equation} 
For $\bar P_{\bar\Omega}(k)$, we estimate
\begin{align}
  \bar P_{\bar\Omega}(k) &\le \sum_{n=2}^\infty\frac{|k|^{ns}}{n!}
  \int_{\mathbb{R}^d\setminus\Omega}\left|f(x)\right|^n\D x \nonumber\\ &\le
  \sum_{n=2}^\infty\frac{|k|^{ns}}{n!}\left[
    |f(0)|-\sigma(\varepsilon)
  \right]^{n-2}\underbrace{\int_{\mathbb{R}^d}|f(x)|^2\D x}_{=:\,C\,<\,\infty}
  \nonumber\\ &=
  \frac{C}{\left[
    |f(0)|-\sigma(\varepsilon)
  \right]^2}\left\{
    \E^{-|k|^s\left[f(0)+\sigma(\varepsilon)\right]}-1+|k|^s\left[
      f(0)+\sigma(\varepsilon)
    \right]
  \right\} \nonumber\\ &\sim
  \frac{C}{\left[
    |f(0)|-\sigma(\varepsilon)
  \right]^2}\E^{-|k|^s\left[f(0)+\sigma(\varepsilon)\right]}
  \quad\mbox{for}\quad |k|\to\infty\;.
\label{eq:64}
\end{align}
Again, this result is exponentially suppressed by the factor $\exp(-|k|^s\sigma(\varepsilon))$ compared to the asymptotic expansion (\ref{eq:57}). Since $\tilde f\in L_2(\mathbb{R}^d)$, we can conclude
\begin{align}
  P(k) &= \bar P(k)-|k|^s\tilde f(k) \sim \bar P(k) \sim \bar P_\Omega(k)
  \nonumber\\ &\sim
  \E^{-|k|^sf(0)}\sqrt{\frac{(2\pi)^d}{|k|^{sd}\det A}}
  \exp\left(-\frac{k^\top A^{-1}k}{2|k|^s}\right)
  \quad\mbox{for}\quad |k|\to\infty\;,
\label{eq:65}
\end{align}
as claimed.

\section{Summary}

We have shown here that Laplace-Fourier-type integrals of the form
\begin{equation}
  P(k) = \int_\Omega\E^{-|k|^sf(x)}\E^{\I kx}\,\D x
\label{eq:66}
\end{equation}
behave asymptotically like
\begin{equation}
  P(k) \sim \E^{-|k|^sf(0)}\sqrt{\frac{(2\pi)^d}{|k|^{sd}\det A}}
  \exp\left(-\frac{k^\top A^{-1}k}{2|k|^s}\right)
\label{eq:67}
\end{equation}
for $|k|\to\infty$ under the following conditions:
\begin{enumerate}
  \item The function $f:\Omega\to\mathbb{R}$ has a negative minimum at $x_0 = 0$;
  \item the Hessian $A$ of $f$ in $x_0 = 0$ exists and is positive definite;
  \item $f$ is square integrable on $\Omega$; and
  \item the exponent $s\ge2$.
\end{enumerate}
This result is valid for $\Omega\subseteq\mathbb{R}^d$ with $d\ge1$. We began in Sect.~\ref{sec:2} by applying Erdélyi's theorem \ref{the:1} to the more general class (\ref{eq:15}) of integrals in one dimension, discussed special cases, and restricted the general result to the case (\ref{eq:66}) with $s\ge2$. We extended the discussion to $d>1$ dimensions in Sect.~\ref{sec:3} and finally showed that the integration domain can be $\mathbb{R}^d$.

Integrals of this type with rapidly oscillating integrands frequently occur in physics, most notably in statistical physics. Asymptotic expressions like (\ref{eq:67}) can help to understand the behaviour of such integrals in the limit of small scales or large momenta, $|k|\to\infty$, and as a test case for numerical integrations. In a forthcoming paper, we shall apply our result to the formation of small-scale cosmic structures.

\section*{Acknowledgements}

We thank Manfred Salmhofer for numerous helpful discussions and for reviewing a draft of this paper. Constructive comments by Robert Lilow and Ricardo Waibel contributed to clarifying the line of reasoning. This work is funded by the Deutsche Forschungsgemeinschaft (DFG, German Research Foundation) under Germany's Excellence Strategy EXC 2181/1 - 390900948 (the Heidelberg STRUCTURES Excellence Cluster).

\appendix

\section{Explicit expressions for the $I_m(\lambda, k)$ symbols}

We start with the expression
\begin{equation}
  I_{m}(\lambda, k) = \frac{1}{m! (a_0 \lambda)^{(m+\beta)/\alpha}}
  \sum_{n=0}^{\infty}\frac{\Gamma\left(\frac{n+m+\beta}{\alpha}\right)}{n!} \left(
    \frac{\I k}{(a_0 \lambda)^{1/\alpha}}
  \right)^n d_{m,n+m}
\label{eq:68}
\end{equation}
for these symbols, obtained from (\ref{eq:18}) inserting the definition $\nu = (n+\beta)/\alpha$, containing the coefficients $d_{m,n+m}$ from (\ref{eq:9}) with $d_{0,n} = 1$. For $m>0$, we make use of the result obtained by \cite{belkic2019} for the expansion of $\D^m/ \D x^m \left[1/f^{\lambda}(x)\right]$,
\begin{align}
  d_{m,n+m} &= \frac{\D^m}{\D x^m} \left(
    \frac{1}{1+\sum_{j=1}^\infty \frac{a_j}{a_0} x^j}
  \right)^{-\frac{n+m+\beta}{\alpha}}_{x=0} \nonumber \\ &=
  \sum_{j=1}^m(-1)^j
  \frac{\Gamma\left(\frac{n+m+\beta}{\alpha}+j\right)}
       {\Gamma\left(\frac{n+m+\beta}{\alpha}\right)}
  B_{m,j}\left(\frac{a_1}{a_0},\dotsm,\frac{a_{m-j+1}}{a_0}\right)\;.
\label{eq:69}
\end{align}
Here, $B_{m,j}$ denote the exponential Bell polynomials
\begin{equation}
    B_{m,j}(x_1,\dotsc,x_{m-j+1}) = \sum \frac{m!}{j_1!j_2!\dotsm j_{m-j+1}! } x_1^{j_1} x_{2}^{j_2} \dotsm x_{m-j+1}^{j_{m-j+1}}\;,
    \label{eq:70}
\end{equation}
where the summation extends over the non-negative integers $j_1,j_2,\dotsc,j_{m-j+1}$ subject to the conditions
\begin{align}
  j_1+j_2+\dotsm+j_{m-j+1}&=j\;, \nonumber \\
  j_1+2 j_2+\dotsm + (m-j+1)j_{m-j+1} &= m\;.
\label{eq:71}
\end{align}
Combining these results, we arrive at
\begin{align}
  I_{m}(\lambda, k) &= \frac{1}{m!(a_0\lambda)^{(m+\beta)/\alpha}}
  \sum_{n=0}^{\infty} \sum_{j=1}^m(-1)^j
  \frac{\Gamma\left(\frac{n+m+\beta}{\alpha}+j\right)}{n!}
  \left(\frac{\I k}{(a_0\lambda)^{1/\alpha}}\right)^n
  B_{m,j}\left(\frac{a_1}{a_0},\dotsm,\frac{a_{m-j+1}}{a_0}\right)
  \nonumber \\ &=
  \frac{1}{m!(a_0\lambda)^{(m+\beta)/\alpha}}
  \sum_{j=1}^m\left(-\frac{1}{a_0}\right)^jB_{m,j}
  \left(a_1,\dotsm,a_{m-j+1}\right)
  \sum_{n=0}^{\infty}\frac{\Gamma\left(\frac{n+m+\beta}{\alpha}+j\right)}{n!}
  \left(\frac{\I k}{(a_0\lambda)^{1/\alpha}}\right)^n \nonumber \\ &=
  \frac{1}{m!(a_0\lambda)^{(m+\beta)/\alpha}}
  \sum_{j=1}^m\left(-\frac{1}{a_0}\right)^j B_{m,j}
  \left(a_1,\dotsm,a_{m-j+1}\right)
  {}_1\Psi_0\left(
    \begin{smallmatrix}
      \left((m+\beta)/\alpha+j, 1/\alpha\right) \\ -
    \end{smallmatrix};
    \frac{\I k \mu}{(a_0 \lambda)^{1/\alpha}}
  \right)\;,
\label{eq:72}
\end{align}
where ${}_1\Psi_0$ denotes the Fox-Wright function \cite{fox1928, wright1935, wright1940}.

\bibliography{main}

\end{document}